\documentclass[12pt]{article}
\usepackage[cp866]{inputenc}
\usepackage{textcomp}
\usepackage{amsfonts,amssymb,amsmath,fullpage,pst-node,rotating,graphics,epsfig,subfigure,amsthm}
\usepackage{multirow}
\usepackage{authblk}

\newcommand{\ri}{{\rm i}\,}

\newcommand{\re}{{\rm e}}
\newcommand{\rd}{\,{\rm d}\,}

\newcommand{\bff}{{\bf f}}
\newcommand{\bF}{{\bf F}}

\newcommand{\bu}{{\bf u}}
\newcommand{\bT}{{\bf T}}
\newcommand{\bP}{{\bf P}}
\newcommand{\bM}{{\bf M}}
\newcommand{\be}{{\bf e}}

\newcommand{\bx}{{\bf x}}
\newcommand{\bv}{{\bf v}}
\newcommand{\bb}{{\bf b}}
\newcommand{\bn}{{\bf n}}

\newcommand{\bd}{{\bf d}}

\newcommand{\cH}{{\cal H}}
\newcommand{\cA}{{\cal A}}
\newcommand{\cB}{{\cal B}}
\newcommand{\cC}{{\cal C}}
\newcommand{\rme}{{\rm e}}
\newcommand{\rr}{{\rm r}}

\newcommand{\cW}{{\cal W}}
\newcommand{\cU}{{\cal U}}
\newcommand{\cV}{{\cal V}}

\newtheorem{theorem}{Theorem}

\newtheorem{remark}{Remark}
\newtheorem{lemma}{Lemma}

\begin{document}

\title{An efficient Galerkin method for problems with physically
realistic boundary conditions.}

\author{Olga Podvigina\\
Institute of Earthquake Prediction Theory\\
and Mathematical Geophysics\\
84/32 Profsoyuznaya St., 117997 Moscow, Russian Federation}

\maketitle

\begin{abstract}
The Galerkin method is often employed for numerical integration of evolutionary
equations, such as the Navier--Stokes equation or the magnetic induction
equation. Application of the method requires solving an equation of the form
$P(Av-f)=0$ at each time step, where $v$ is an element of a
finite-dimensional space $\cV$ with a basis satisfying boundary
conditions, $P$ is the orthogonal projection on this space and $A$ is a linear
operator. Usually the coefficients of $v$ expanded in the basis are found by
calculating the matrix of $PA$ acting on $\cV$ and solving the respective
system of linear equations. For physically realistic boundary conditions
(such as the no-slip boundary conditions for the velocity, or for a dielectric
outside the fluid volume for the magnetic field) the basis is often
not orthogonal and solving the problem can be computationally demanding.
We propose an algorithm giving an opportunity to
reduce the computational cost for such a problem. Suppose there exists
a space $\cW$ that contains $\cV$, the difference between the dimensions of
$\cW$ and $\cV$ is small relative to the dimension of $\cV$, and solving
the problem $P(Aw-f)=0$, where $w$ is an element of $W$, requires less
operations than solving the original problem. The equation $P(Av-f)=0$
is then solved in two steps: we solve the problem $P(Aw-f)=0$ in $\cW$, find
a correction $h=v-w$ that belongs to a complement to $\cV$ in $\cW$, and
obtain the solution $w+h$.
When the dimension of the complement is small the proposed algorithm is more
efficient than the traditional one.
\end{abstract}

\section{Introduction}
Spectral methods are often employed for numerical solution of partial
differential equations, because relying on a relatively small number
of unknowns they provide accurate approximations to solutions \cite{boy,can,got,she11}.
The three most commonly used spectral methods are the Galerkin, tau and
collocation methods. In the Galerkin and tau methods the solution is expanded
in series of basis functions. The coefficient of the series are found from the condition
that the scalar product of  the residual with
each test function is zero. In the collocation method the solution is
approximated in the physical space by minimising the discrepancy over
a set of collocation points.

Solving an evolutionary equation
$${\partial\over\partial t}\bv(\bx,t)=F(\bv(\bx,t))$$
by the Galerkin method requires finding the
coefficients $\{c_j(t)\}_{j=1}^N$ such that the approximate solution
$$\bv(\bx,t)=\sum_{j=1}^Nc_j(t)\phi_j(\bx)$$
satisfies
\begin{equation}\label{gal0}
({\partial\over\partial t}\bv(\bx,t)-F(\bv(\bx,t)),\phi_j(\bx))
\hbox{ for any }1\le j\le N\hbox{ and }0<t<T.
\end{equation}
Here the functions $\{\phi_j(\bx)\}_{j=1}^N$ satisfy the desired boundary
conditions, and $(\cdot,\cdot)$ is a scalar product in the Hilbert space
$\cH$ to which they belong.

Denote by $\cV$ the space spanned by $\{\phi_j(\bx)\}_{j=1}^N$ and
by $P_{\cV}$ the orthogonal projection on $\cV$.
Provided $\bv(\bx,0)$ belongs to $\cV$, equation (\ref{gal0}) is equivalent to
$${\partial\over\partial t}\bv(\bx,t)=P_{\cV}F(\bv(\bx,t))\hbox{ for any }
0<t<T.$$

Depending on the equation and boundary conditions, the basis
$\{\phi_j(\bx)\}_{j=1}^N$ can be comprised of Fourier exponentials, linear
combinations of Chebychev, Legendre, Jacobi or some other polynomials,
or their products \cite{can,haf20,guoa,guob,she94,spa91,tem95}.
The projection $P_{\cV}$ can be easily computed for the conditions of space
periodicity, since the Fourier exponentials comprise an orthogonal basis and they
individually satisfy the conditions. For other boundary
conditions, the functions $\phi_j$ usually
involve linear combinations of polynomials and computing
$P_{\cV}F(\bv(\bx,t))$ may become a non-trivial problem, which is extensively
studied in literature \cite{can,leo82,mos83,spa86,zhe}.

In this paper we propose a new algorithm for computing
$P_{\cV}F(\bv(\bx,t))$ or, in more general terms, for computing $v$ such that
\begin{equation}\label{gal2}
P_{\cV}(\cA v-f)=0\hbox{ and }v\in\cV,
\end{equation}
where $\cA$ is a linear operator.

The algorithm resembles the tau method: First, we compute $w$,
an approximation to $v$, such that
\begin{equation}\label{gal3}
P_{\cV}(\cA w-f)=0\hbox{ and }w\in\cW,
\end{equation}
where $w$ is not required to satisfy boundary conditions.
Second, we add to $w$ a correction $h$ such that $v=w+h$ solves (\ref{gal2}).
The space $\cW\supset\cV$ must satisfy two conditions: (i) it is more
convinient to compute $w$ than to compute $v$, and (ii) the dimension of
$\cW\ominus\cV$ is small.
This is still (a varient of) the Galerkin method, because the same solution is
obtained, as by the traditional Galerkin algorithm (up to roundoff errors).

In the next section we prove a theorem providing a mathematical background
for the algorithm and define the algorithm as a sequence of steps.
In Section \ref{sec3} we apply the algorithm for solving equations (\ref{gal2})
in the context of one-dimensional problems stated on the interval $[-1,1]$
subject to the boundary conditions $v(\pm1)=0$; the basis in $\cV$ is then
comprised of linear combinations of Chebyshev polynomials.

In section \ref{sec4} we
apply the algorithm for solving equations of magnetoconvection
in a plane horizontal layer. The flow velocity, temperature and magnetic field
are required to be periodic in horizontal directions. The fields are
decomposed in series in the Fourier exponentials in $x_1$ and $x_2$
and Chebyshev polynomials in $x_3$. Computing $P_{\cV}F(\bv(\bx,t))$
is reduced to solving a number of one-dimensional problems, some of which
have been solved in section \ref{sec3}.

In the conclusion we briefly summarise the results and discuss other
possible implementations of the algorithm.

\section{Description of the algorithm}\label{sec2}

Recall the traditional formulation of the Galerkin method.
Let $\cH$ be a Hilbert space with a basis $\{\phi_i\}^{\infty}_{i=1}$, $\cA$ an
operator acting on $\cH$ and $f$ an element of $\cH$. Solving the problem
$\cA v=f$ by the Galerkin method requires finding a set of coefficients
$\{c_i\}^N_{i=1}$ such that the sum
$$v=\sum_{i=1}^Nc_i\phi_i$$
satisfies the relations
$$(\cA v-f,\phi_i)=0,\hbox{ for any }1\le i\le N.$$
Here $(\cdot,\cdot)$ denotes the scalar product in $\cH$.

The method can be formulated in a slightly different way. Denote by $\cV$ the
subspace of $\cH$ spanned by $\{\phi_i\}^N_{i=1}$ and
by $P_{\cV}$ the orthogonal projection on $\cV$. Solving the problem $\cA v=f$
by the Galerkin method requires finding $v$ such that
\begin{equation}\label{metG}
P_{\cV}(\cA v-f)=0\hbox{ and }v\in\cV.
\end{equation}
Evidently, the two formulations are equivalent.

Below we prove a theorem that gives the mathematical foundation of the
new algorithm. The proof relies on the following lemma:

\begin{lemma}\label{lem1}
Let $\cB$ be a linear operator acting on a finite-dimensional space $\cW$,
$\cV$ a subspace of $\cW$ and $P_{\cV}$ an orthogonal projection on $\cV$.
Denote by $\cV^{\perp}$ the orthogonal complement to $\cV$ in $\cW$.
Suppose the operator $\cC=P_{\cV}\cB:\,\cV\to \cV$ satisfies $\det\cC\ne0$.
(Here we use the same notation for the linear operator $\cC$ and for its
matrix.) Then any $s\in \cV^{\perp}$ can be uniquely decomposed as
\begin{equation}\label{decos}
s=q+\tilde q,\hbox{ where }P_{\cV}\cB q=0
\hbox{ and }\tilde q\in\cV.
\end{equation}
\end{lemma}

\proof
Denote $f=P_{\cV}\cB s$. If $f=0$, then
in the decomposition (\ref{decos}) we set $q=s$ and $\tilde q=0$.

If $f\ne0$, we define $\tilde q\in{\cV}$ as the solution to the equation
$\cC\tilde q=f$. The solution exists and is unique due to the condition that
$\det\cC=0$. The difference $q=s-\tilde q$ satisfies
$$P_{\cV}\cB q=P_{\cV}\cB(s-\tilde q)=P_{\cV}\cB s-P_{\cV}\cB\tilde q=
f-\cC\tilde q.$$
Therefore, $P_{\cV}\cB q=0$.
\qed

\begin{theorem}\label{th0}
Let the operator $\cB$ and the spaces $\cW$ and $\cV$ be the same as in lemma
\ref{lem1} and $f$ be an element of $\cW$. Denote by $s_i$, $1\le i\le K$,
vectors comprising an orthonormal basis in $\cV^{\perp}$. By lemma \ref{lem1},
$s_i=q_i+\tilde q_i$, where $P_{\cV}\cB q_i=0$ and $\tilde q_i\in\cV$.
Suppose that $w\in\cW$ is a solution
\footnote{Note, that the solution is not unique.}
to the problem
$$P_{\cV}(\cB w-f)=0.$$
Then
\begin{equation}\label{corrstep}
v=w-\sum_{i=1}^Kb_iq_i,\hbox{ where }b_i=(w,s_i),\  i=1,...,K,
\end{equation}
satisfies
$$P_{\cV}(\cB v-f)=0\hbox{ and }v\in{\cV}.$$
\end{theorem}

\proof
The inclusion $v\in{\cV}$ holds true because
$$v=w-\sum_{i=1}^Kb_iq_i=w-\sum_{i=1}^Kb_i(s_i-\tilde q_i)=
(w-\sum_{i=1}^K(w,s_i)s_i)+\sum_{i=1}^Kb_i\tilde q_i,$$
where the first term in the sum belongs to ${\cV}$ since the basis $s_i$ is
an orthonormal basis in ${\cV}^{\perp}$, and the second term belongs to ${\cV}$
since all $\tilde q_i$ are elements of ${\cV}$. The equation
$P_{\cV}(\cB v-f)=0$ is satisfied, because
$$P_{\cV}(\cB v-f)=P_{\cV}(B(w-\sum_{i=1}^Kb_iq_i)-f)=P_{\cV}(Bw-f)-
\sum_{i=1}^Kb_iP_{\cV}(\cB q_i)=0.$$
\qed

The new algorithm is applicable when $\cA$ in (\ref{metG}) is a linear
operator. Suppose that there exists a space $\cW\supset{\cV}$, such that
solving the problem $P_{\cV}(\cA w-f)=0$, where $w\in \cW$, is
more efficient than solving $P_{\cV}(\cA v-f)=0$, where $v\in{\cV}$.
The proposed algorithm is comprised of three steps.\\
{\it Preliminary step.} Computing an orthonormal basis in ${\cV}^{\perp}$
and decomposing the basis as $s_i=q_i+\tilde q_i$, following
\footnote
{The lemma requires that $\det\cC\ne0$, where $\cC$ is the matrix of the operator
$P_{\cV}\cA:{\cV}\to{\cV}$. Generically the condition is satisfied.
This is always implicitly assumed when using the traditional Galerkin
method, because otherwise the equation (\ref{metG}) does not have
a solution for certain $f$'s.}
lemma \ref{lem1}. The basis is found by solving the equation $P_{\cV}s_i=0$,
where $s_i\in \cW$, and the vectors $\tilde q$ by solving
$P_{\cV}\cA(\tilde q_i-s_i)=0$, where $\tilde q_i\in{\cV}$.\\
{\it Main step.} Solving the equation $P_{\cV}(\cA w-f)=0$, where $w\in \cW$.\\
{\it Correction step.} Computing $v$ from $w$ by (\ref{corrstep}).

The method is efficient if the dimension of ${\cV}^{\perp}$ is considerably
smaller than that of ${\cV}$ and equation (\ref{metG}) is solved many times
for different $f$'s. This happens, e.g., in the course of numerical integration
of evolutionary equations. An example of such an application
is presented in section \ref{sec4}. The preliminary step is performed only
once and the obtained vectors $s_i$ and $q_i$ are stored. Since
$\dim {\cV}^{\perp}\ll\dim {\cV}$, the memory usage increases insignificantly
and the correction step requires a small amount of computations.
The efficiency is achieved because the problem
$P_{\cV}(\cA w-f)=0$, where $w\in\cW$, is easier to solve than the original one.

\begin{remark}\label{rem1}
Solving the problem discussed in the beginning of this section by
the Petrov--Galerkin method requires finding $v$ such that
$$P_{\cU}(\cA v-f)=0\hbox{ and }v\in{\cV}.$$
Here $\cU$ is a subspace of $\cH$, different from $\cV$, and $\dim\cU=\dim\cV$.
In other words, the only difference between the Galerkin and the Petrov--Galerkin
method is that $P_{\cV}$ is replaced by $P_{\cU}$. The lemma and the theorem
hold true for $P_{\cV}$ replaced by $P_{\cU}$. Thus, the algorithm remains
applicable for solving linear equations by the Petrov--Galerkin method.
\end{remark}

\section{Two simple examples}\label{sec3}

\subsection{Example 1}\label{sec31}

{\it The problem.}
Let $\cH$ be the space of $L^2_w$ functions on the interval $[-1,1]$ with
the basis
comprised of Chebyshev polynomials. For the scalar product defined as
$$(f,g)=\int_{-1}^1{fg\over\sqrt {1-x^2}}\rd x,$$
the basis is orthogonal. We choose the space $\cW$ spanned
by a finite number of polynomials, $\cW=<T_0,T_1,...,T_{N+1}>$, and $\cV$
is the subspace of $\cW$ comprised of polynomials vanishing at $\pm1$.
A basis in $\cV$ is, e.g.,
$T_2-T_0,T_3-T_1,...,T_{N+1}-T_1$. (Here and below an even $N$ is assumed.)

The operator $\cA$ is the identity. The equation (\ref{metG}) takes now the form
\begin{equation}\label{metGI}
P_{\cV}(v-f)=0\hbox{ and }v\in{\cV}.
\end{equation}

\noindent
{\it The preliminary step.}
Since $\cA$ is the identity, $P_{\cV}\cA{\cV}^{\perp}=0$, and thus
$\tilde q_i=0$. Since
\begin{equation}\label{prodT}
(T_n,T_m)=\left\{
\begin{array}{l}
0\quad n\ne m\\
\pi\quad n=m\ne0\\
\pi/2\quad n=m=0
\end{array}
\right.,
\end{equation}
the basis in ${\cV}^{\perp}$ comprised of
\begin{equation}\label{basP}
\renewcommand{\arraystretch}{1.8}
\begin{array}{l}
s_1=q_1=\displaystyle{4\over\sqrt{\pi}(N+4)}(2T_0+T_2+...+T_N),\\
s_2=q_2=\displaystyle{4\over\sqrt{\pi}(N+2)}(T_1+T_3+...+T_{N+1})
\end{array}
\end{equation}
is orthonormal.

\noindent
{\it The main step.}
For a given
$$f=\sum_{i=0}^M f_iT_i,$$
a solution to $P_{\cV}(w-f)=0$, where $w\in \cW$, is
$$w=\sum_{i=0}^{N+1} f_iT_i.$$

\noindent
{\it The correction step.}
By (\ref{prodT}) and (\ref{basP}),
$$
b_1=(w,s_1)=\displaystyle{2\sqrt{\pi}\over N+4}
(2f_0+\sum_{i=2,\ i\hbox{ even}}^N f_i),\quad
b_2=(w,s_2)=\displaystyle{2\sqrt{\pi}\over N+2}
\sum_{i=1,\ i\hbox{ odd}}^{N+1} f_i.$$
Finally,
$$v=w-b_1s_1-b_2s_2.$$

Hence, computing $v$ is reduced to computing two scalar products and two
subtractions.
By contrast, application of the traditional algorithm requires
solving a system of $N$ linear equations in $N$ variables \cite{zhe}.

\subsection{Example 2}\label{sec32}

{\it The problem.} The spaces $\cH$, $\cW$ and $\cV$ are the same, as
in example 1.
We solve equation ({\ref{metG}) for $\cA v=\alpha v+\beta v''$.

\noindent
{\it The preliminary step.}
The basis (\ref{basP}) in ${\cV}^{\perp}$ is used. The vectors
$\tilde q_i\in{\cV}$ are solutions to the problem $P_{\cV}\cA(\tilde q_i-s_i)=0$
obtained by any standard algorithm.
The vectors $q_1$ and $q_2$ satisfy $q_j=s_j-\tilde q_j$.

\noindent
{\it The main step.}
The element $f$ and the solution $w\in W$ are expanded in sums
of Chebyshev polynomials:
$$f=\sum_{i=0}^M f_iT_i,\quad w=\sum_{i=0}^{N+1} w_iT_i.$$
The equation $P_{\cV}(\cA w-f)=0$ takes the form
\begin{equation}\label{eqsum}
\alpha\sum_{i=0}^{N+1} w_iT_i+\beta\sum_{i=0}^{N+1}w_iT''_i=\sum_{i=0}^{N+1}
f_iT_i.
\end{equation}
To set up a recursive algorithm, we modify the equation by adding an extra
term,
\begin{equation}\label{eqsumn}
\alpha\sum_{i=0}^{N+1} w_iT_i+\beta\sum_{i=0}^{N+1}w_iT''_i=
\sum_{i=0}^{N+1} f_iT_i+\sum_{i=0}^{N+1} a_iT''_i,
\end{equation}
where initially we set $a_i=0$ for all $0\le i\le N+1$.
Since $T''_i$ is a sum of polynomials of degree $i-2$ or less, (\ref{eqsum})
implies
$$w_{N+1}={1\over\alpha}f_{N+1}\hbox{ and }w_N={1\over\alpha}f_N.$$
Using the relation
$$2T_k={1\over2(k+1)(k+2)}T''_{k+2}-{1\over(k-1)(k+1)}T''_k+
{1\over2(k-1)(k-2)}T''_{k-2},$$
we replace $T''_{N+1}$ and $T''_N$ in (\ref{eqsumn}) by sums of polynomials
and derivatives of smaller indices. In the modified equations,
$$f_{N-1}^{new}=f_{N-1}^{old}-4\beta N(N+1)w_{N+1},\
f_{N-2}^{new}=f_{N-2}^{old}-4\beta N(N-1)w_N,$$
$$a_{N-1}^{new}=a_{N-1}^{old}+{2(N+1)\over N-2}(\beta w_{N+1}-a_{N+1}),\
a_{N-2}^{new}=a_{N-2}^{old}+{2N\over N-3}(\beta w_N-a_N),$$
$$a_{N-3}^{new}=a_{N-3}^{old}-{N(N+1)\over (N-2)(N-3)}(\beta w_{N+1}+a_{N+1}),\
a_{N-4}^{new}=a_{N-4}^{old}-{N(N-1)\over (N-3)(N-4)}(\beta w_N+a_N),$$
while no other terms are changed.

As a result of this substitution, we obtained the same equation as
(\ref{eqsumn}), except that summation is performed up to $i=N-1$.
Repeating the above calculations with indices reduced by two, we
obtain $w_{N-1}$, $w_{N-2}$ and a new equation with summation up to $N-3$.
Repeating this procedure $N/2-1$ more times, we calculate all coefficients $w_i$.

\noindent
{\it The correction step} is performed as described in theorem \ref{th0}.

\section{Application for solving the equations of
magnetohydrodynamics}\label{sec4}

\subsection{Equations}\label{sec41}

We consider a conductive fluid in a plane horizontal layer heated from
below and rotating about an inclined axis. The equations governing the temporal
evolution of the fluid and the magnetic field are the Navier--Stokes equation
\begin{equation}\label{nst}
{\partial{\bf v}\over\partial t}={\bf v}\times(\nabla\times{\bf v})
+P\tau{\bf v}\times{\bf e}_{\rr}
+P\Delta{\bf v}+PR\theta{\bf e}_z-\nabla p-{\bf b}\times(\nabla\times{\bf b}),
\end{equation}
the magnetic induction equation
\begin{equation}\label{mind}
{\partial{\bf b}\over\partial t}=\nabla\times({\bf v}\times{\bf b})
+\eta\Delta{\bf b}
\end{equation}
and the heat transfer equation
\begin{equation}\label{htr}
{\partial\theta\over\partial t}=-({\bf v}\cdot\nabla)\theta+v_z+\Delta\theta,
\end{equation}
where ${\bf e}_{\rr}$ is the unit vector parallel to the rotation axis.
The flow and the magnetic field are solenoidal
\begin{equation}\label{inc0}
\nabla\cdot{\bf v}=\nabla\cdot{\bf b}=0.
\end{equation}
Here $\bf v$ denotes the flow velocity, $p$ the modified pressure, $\bf b$
the magnetic field and $\theta$ the difference between the flow temperature
and the linear temperature profile.
The equations involve the following parameters: the Prandtl number, $P$, and
the magnetic Prandtl number, $Pm$, describing the properties of
the fluid, the Rayleigh number, $R$, proportional to the temperature difference
between the upper and lower boundaries, and the Taylor number, $Ta=\tau^2$,
where $\tau$ is proportional to angular velocity.

The flow satisfies the no-slip condition on the horizontal boundaries,
which are kept at constant temperatures:
\begin{equation}\label{bou1}
v_1=v_2=v_3=0,\quad\theta=0\quad\hbox{ at }x_3=\pm1.
\end{equation}
We assume electrically insulating upper and perfectly electrically conducting
lower boundaries. The boundary conditions for the magnetic field are, thus,
\begin{equation}\label{boum}
\renewcommand{\arraystretch}{1.8}
\begin{array}{l}
{\bf b}|_{x_3=1}=\nabla h,\quad
\nabla^2 h=0\hbox{ for }x_3\ge 1,\quad
h\to0\hbox{ for }x_3\to\infty,\\
\displaystyle{{\partial b_1\over\partial x_3}|_{x_3=-1}=
{\partial b_2\over\partial x_3}|_{x_3=-1}=
b_3|_{x_3=-1}=0.}
\end{array}
\end{equation}
The fields are periodic in the horizontal directions:
\begin{equation}\label{hor}
\renewcommand{\arraystretch}{1.3}
\begin{array}{l}
{\bf v}(\bx)=\bv(x_1+n_1L_1,x_2+n_2L_2,x_3),\\
{\bf b}(\bx)=\bb(x_1+n_1L_1,x_2+n_2L_2,x_3),\\
\theta(\bx)=\theta(x_1+n_1L_1,x_2+n_2L_2,x_3).
\end{array}
\end{equation}

\subsection{An explicit method}\label{sec42}

In this section we give an example of application of the new algorithm
for numerical integration of equations (\ref{nst})--(\ref{htr}) by an
explicit time stepping. For simplicity, the Euler scheme is used.
Generalisations for more complex methods, such as the Runge--Kutta or
Adams--Bashford methods, can be done with minor modifications.

The system (\ref{nst})--(\ref{htr}) takes the form $\dot\bu=f(\bu)$, where
$\bu=(\bv,\theta,\bb)$ and
$f(\bu)=(f_{\bv}(\bu),f_{\theta}(\bu),f_{\bb}(\bu))$.
The Euler scheme with the time step $\delta t$ results in the relation
\begin{equation}\label{intu0}
\bu(t_{m+1})=\bu(t_m)+\delta tf(\bu(t_m)).
\end{equation}
For the unknown fields expanded in finite sums of products of Fourier
exponentials in $x_1$ and $x_2$ and Chebyshev polynomials in $x_3$,
the Galerkin method stipulates modifying equation (\ref{intu0}) into the form
\begin{equation}\label{intu1}
\bu(t_{m+1})=\bu(t_m)+\delta tP_{\cU}f(\bu(t_m)),
\end{equation}
where $P_{\cU}$ is the orthogonal projection
into space of fields satisfying the boundary conditions.

In the following subsections we perform the spatial
discretisation and derive numerical schemes for temperature, magnetic
field and flow velocity individually.

\subsubsection{The temperature}\label{sec421}

The temperature field is sought as a sum of products of the Fourier
exponentials of $x_1$ and $x_2$ and functions of $x_3$,
\begin{equation}\label{temp1}
\theta(\bx,t)=\sum_{|n_1|\le N_1,|n_2|\le N_2}
\Theta_{n_1,n_2}(x_3,t)\rme^{\ri(\alpha_1n_1x_1+\alpha_2n_2x_2)},
\end{equation}
where each $\Theta_{n_1,n_2}(x_3,t)$ is a sum of Chebyshev polynomials,
\begin{equation}\label{temp2}
\Theta_{n_1,n_2}(x_3,t)=\sum_{0\le n_3\le N_3+1}\theta_{\bn}(t)T_{n_3}(x_3).
\end{equation}
The boundary condition $\theta(x_1,x_2,\pm1,t)=0$ implies
\begin{equation}\label{bcT}
\Theta_{n_1,n_2}(\pm1,t)=0\hbox{ for all }n_1\hbox{ and }n_2.
\end{equation}

For the $\theta$ component of $\bu$, equation (\ref{intu1}) takes the form
\begin{equation}\label{intt}
\theta(t_{m+1})=\theta(t_m)+\delta t P_{\cU}f_{\theta}(\bu(t_m)),
\end{equation}
where $f_{\theta}$ is the r.h.s. of (\ref{htr}) and $P_{\cU}$ is the projection
into the space $\cU$ of functions (\ref{temp1}),(\ref{temp2}) satisfying
boundary conditions (\ref{bcT}).

The nonlinear term in (\ref{htr}) is computed in the physical space
employing an FFT to switch between the physical and spectral spaces.
Differentiation in $x_1$ and $x_2$ is trivial. To differentiate in $x_3$,
we use the relation
$$2T_k(z)={1\over k+1}T'_{k+1}(z)-{1\over k-1}T'_{k-1}(z)$$
implying that the coefficients of the series
\begin{equation}\label{difz1}
a(z)=\sum_{0\le k\le K+1}a_kT_k(z)\hbox{ and }
a'(z)=b(z)=\sum_{0\le k\le K}b_kT_k(z)
\end{equation}
satisfy the relation \cite{can}
\begin{equation}\label{difz2}
b_k=b_{k+2}+2(k+1)a_{k+1},\ 1\le k\le K,\quad 2b_0=b_2+2a_1.
\end{equation}
Therefore, the coefficients $b_k$ can be computed for decreasing $k$
starting from $k=K$.

We expand $f_{\theta}$ similarly to $\theta$,
\begin{equation}\label{ff1}
f_{\theta}(\bx,t)=\sum_{|n_1|\le K_1,|n_2|\le K_2}
F_{n_1,n_2}(x_3,t)\rme^{\ri(\alpha_1n_1x_1+\alpha_2n_2x_2)},
\end{equation}
where
\begin{equation}\label{ff2}
F_{n_1,n_2}(x_3,t)=\sum_{0\le n_3\le K_3+1}f_{\bn}(t)T_{n_3}(x_3).
\end{equation}
By (\ref{bcT}), when computing $P_{\cU}f_{\theta}$, the projection $P_{\cU}$
into the space of fields satisfying the boundary condition can be computed
individually for each term $F_{n_1,n_2}(x_3,t)$. The numerical scheme
(\ref{intt}) now takes the form
\begin{equation}\label{intt12}
\Theta_{n_1,n_2}(x_3,t_{m+1})=\Theta_{n_1,n_2}(x_3,t_m)+
\delta t P_{\cV}(F_{n_1,n_2}(x_3,t_m))
\end{equation}
for all $(n_1,n_2)$, where $\cV$ and $P_{\cV}$ are the same as in example 1.

\subsubsection{The magnetic field}\label{sec422}

The magnetic field is represented as the sum of the toroidal, $\bT^b$,
poloidal, $\bP^b$, and mean field, $\bM^b$, components,
\begin{equation}\label{magf0}
\bb(\bx,t)=\bT^b(\bx,t)+\bP^b(\bx,t)+\bM^b(x_3,t),
\end{equation}
where
\begin{equation}\label{magf1}
\renewcommand{\arraystretch}{1.5}
\begin{array}{l}
\bT^b(\bx,t)=\nabla\times(T^b(\bx,t)\be_3),\\
\bP^b(\bx,t)=\nabla\times\nabla\times(P^b(\bx,t)\be_3),\\
\bM^b(x_3,t)=(M_1^b(x_3,t),M_2^b(x_3,t),0)
\end{array}
\end{equation}

Boundary conditions (\ref{boum}) for the magnetic field imply that the
component satisfy the following boundary
conditions
\begin{equation}\label{magBC}
\renewcommand{\arraystretch}{1.5}
\begin{array}{l}
T^b|_{x_3=1}=0,\ \partial T^b/\partial x_3|_{x_3=-1}=0,\\
(\partial^2 P^b/\partial x_1\partial x_3,
\partial^2 P^b/\partial x_2\partial x_3,
-\partial^2 P^b/\partial x_1^2-\partial^2 P^b/\partial x_2^2)|_{x_3=1}=\nabla h,\\
\partial^2 P^b/\partial x_3^2|_{x_3=-1}=0,\ P^b|_{x_3=-1}=0,\\
\bM^b|_{x_3=1}=0,\   \partial \bM^b/\partial x_3|_{x_3=-1}=0.
\end{array}
\end{equation}

Denote by $\bff_b(\bx,t)$ the r.h.s. of equation (\ref{mind}). It is
computed the same way as the r.h.s. of the heat transfer equation.
As a result, we obtain $\bff_b(\bx,t)$ in the form of the series
\begin{equation}\label{fdeco}
\bff_b(\bx,t)=
\sum_{|n_1|\le N_1,|n_2|\le N_2}
\bF_{n_1,n_2}(x_3,t)\rme^{\ri(\alpha_1n_1x_1+\alpha_2n_2x_2)},
\end{equation}
where
$\bF_{n_1,n_2}(x_3,t)=(F^1_{n_1,n_2}(x_3,t),F^2_{n_1,n_2}(x_3,t),F^3_{n_1,n_2}(x_3,t))$
and each $F^j_{n_1,n_2}(x_3,t)$ is decomposed as a sum of $N_3+2$
Chebyshev polynomials in $x_3$.

Since $\bff_b(\bx,t)$ is solenoidal, it is a sum of the poloidal, toroidal and
mean-field components:
\begin{equation}\label{magrhs}
\bff_b(\bx,t)=\bT^f(\bx,t)+\bP^f(\bx,t)+\bM^f(x_3,t).
\end{equation}
The components are expressed in terms of $T^f(\bx,t)$, $P^f(\bx,t)$,
$M_1^f(x_3,t)$ and $M_2^f(x_3,t)$ as in (\ref{magf1}).
The toroidal and poloidal components are expanded in the series
\begin{equation}\label{tpdeco}
\renewcommand{\arraystretch}{1.5}
\begin{array}{l}
T^f(\bx,t)=\sum_{|n_1|\le N_1,|n_2|\le N_2}
G^T_{n_1,n_2}(x_3,t)\rme^{\ri(\alpha_1n_1x_1+\alpha_2n_2x_2)},\\
P^f(\bx,t)=\sum_{|n_1|\le N_1,|n_2|\le N_2}
G^P_{n_1,n_2}(x_3,t)\rme^{\ri(\alpha_1n_1x_1+\alpha_2n_2x_2)}.
\end{array}
\end{equation}
By (\ref{magf1}), (\ref{magrhs}) and (\ref{tpdeco}), the coefficients of the
components can be determined from the relations
\begin{equation}\label{detcomp}
\renewcommand{\arraystretch}{1.8}
\begin{array}{l}
-(\alpha_1^2n_1^2+\alpha_2^2n_2^2)G^T_{n_1,n_2}=(\ri \alpha_2n_2 F^1_{n_1,n_2}-\ri \alpha_1n_1 F^2_{n_1,n_2}),\\
-(\alpha_1^2n_1^2+\alpha_2^2n_2^2)G^P_{n_1,n_2}=F^3_{n_1,n_2},\
 M_1=F^1_{0,0},\ M_2=F^2_{0,0}.
\end{array}
\end{equation}

$P_{\cU}\bff_b(\bx,t)$ is computed individually for each component and each
pair $(n_1,n_2)$. We denote by $\cV$ the spaces of $x_3$-dependent functions
$g(x_3)$ such that $g(x_3)\exp(\ri(\alpha_1n_1x_1+\alpha_2n_2x_2))$ satisfies
the boundary conditions (\ref{magBC}). The spaces are different for different
components, and for the poloidal component the spaces depend on $n_1$ and $n_2$.
To simplify the notation, we do not use indices indicating the dependencies.
The space $\cW$ is spanned by the Chebyshev polynomials
of indices from 0 to $N_3+1$.

Since  $T'_j(-1)=(-1)^jj^2$, boundary conditions (\ref{magBC}) imply that
for the toroidal component and any $(n_1,n_2)$ the orthogonal complement
to $\cV$ in $\cW$ is spanned by
\begin{equation}\label{bassT}
\tilde s_1=(1,2,...,2)\hbox{ and }
\tilde s_2=(0,-1,4,...,(-1)^{N+1}(N+1)^2).
\end{equation}
For given $\tilde s_1$ and $\tilde s_2$, an orthonormal basis
in $\cV^{\perp}$ can be easily computed.
For the mean-field component, the projections $P_{\cV}M_1^f(x_3,t)$ and
$P_{\cV}M_2^f(x_3,t)$ are also computed using this basis.

A harmonic function $h$ can be represented as a sum of terms
\begin{equation}\label{decmT}
h_{n_1,n_2}(\bx)=
C\rme^{\ri(\alpha_1n_1x_1+\alpha_2n_2x_2)-(\alpha_1^2n_1^2+\alpha_2^2n_2^2)^{1/2}x_3}.
\end{equation}
Since
$$T_k'(\pm1)=(\pm1)^kk^2\hbox{ and }T_k''(\pm1)=(\pm1)^k{k^4-k^2\over3},$$
the boundary conditions (\ref{magBC}) for the pair $(n_1,n_2)$ take the form
\begin{equation}\label{magBC12}
\renewcommand{\arraystretch}{2.5}
\begin{array}{l}
\sum_{0\le n_3\le N_3+1}((\alpha_1^2n_1^2+\alpha_2^2n_2^2)^{1/2}-n_3^2)b^P_{\bn}=0,\
\sum_{2\le n_3\le N_3+1}(-1)^{n_3}(n_3^4-n_3^2)b^P_{\bn}=0,\\
\sum_{0\le n_3\le N_3+1}(-1)^{n_3}b^P_{\bn}=0,
\end{array}
\end{equation}
where $b^P_{\bn}$ are the coefficients of $P^b$ decomposed in the basis of
Fourier exponentials and Chebyshev polynomials. Therefore, for a given pair
$(n_1,n_2)$, the orthogonal complement to $\cV$ is spanned by three vectors
$\tilde s_1$, $\tilde s_2$ and $\tilde s_3$; the components are
\begin{equation}\label{bassP}
\renewcommand{\arraystretch}{1.8}
\begin{array}{l}
\tilde s_{1,0}=2(\alpha_1^2n_1^2+\alpha_2^2n_2^2)^{1/2},\
\tilde s_{1,j}=((\alpha_1^2n_1^2+\alpha_2^2n_2^2)^{1/2}-j^2)1,\ j=1,...,N_3,\\
\tilde s_{2,j}=(-1)^j(j^4-j^2),\ j=0,...,N_3,\\
\tilde s_{3,0}=2,\ \tilde s^P_{3,j}=(-1)^j,\ j=1,...,N_3.
\end{array}
\end{equation}
An orthonormal basis in $\cV^{\perp}$ can be found by orthogonalising
and normalising these vectors.

Denote by $B_{n_1,n_2}^T(x_3,t)$ and $B_{n_1,n_2}^P(x_3,t)$ the coefficients
of $T^b(\bx,t)$ and $P^b(\bx,t)$, respectively, expanded in the series in
$\exp(\ri(\alpha_1n_1x_1+\alpha_2n_2x_2))$. The Euler time stepping
$$\bb(\bx,t_{m+1})=\bb(\bx,t_m)+\delta t P_{\cU}\bff_b(\bu(\bx,t_m))$$
for the coefficients exploits the relations
\begin{equation}\label{eulb}
\renewcommand{\arraystretch}{1.8}
\begin{array}{l}
B^T_{n_1,n_2}(x_3,t_{m+1})=B^T_{n_1,n_2}(x_3,t_m)+
\delta t P_{\cV}(G^T_{n_1,n_2}(x_3,t_m)),\\
B^P_{n_1,n_2}(x_3,t_{m+1})=B^P_{n_1,n_2}(x_3,t_m)+
\delta t P_{\cV}(G^P_{n_1,n_2}(x_3,t_m)),
\end{array}
\end{equation}
where $G^T_{n_1,n_2}$ and $G^P_{n_1,n_2}$ are obtained using (\ref{detcomp}).
The basis in $\cV^{\perp}$ is defined by (\ref{bassT}) or (\ref{bassP}).
The scheme for the mean-field component is
\begin{equation}\label{eulbm}
M^b_1(x_3,t_{m+1})=M^b_1(x_3,t_m)+\delta t P_{\cV}(F^1_{0,0}(x_3,t_m)),\
M^b_2(x_3,t_{m+1})=M^b_2(x_3,t_m)+\delta t P_{\cV}(F^2_{0,0}(x_3,t_m)),
\end{equation}
where $\cV^{\perp}$ are spanned by (\ref{bassT}). The projection can be computed
as in Example 1.

\subsubsection{The flow velocity}\label{sec423}

Similarly to magnetic field, the velocity is represented as a sum
of the toroidal, poloidal and mean-field components:
\begin{equation}\label{vel0}
\bv(\bx,t)=\bT^v(\bx,t)+\bP^v(\bx,t)+\bM^v(x_3,t),
\end{equation}
where
\begin{equation}\label{vel1}
\renewcommand{\arraystretch}{1.5}
\begin{array}{l}
\bT^v(\bx,t)=\nabla\times(T^v(\bx,t)\be_3),\\
\bP^v(\bx,t)=\nabla\times\nabla\times(P^v(\bx,t)\be_3),\\
\bM^v(x_3,t)=(M_1^v(x_3,t),M_2^v(x_3,t),0)
\end{array}
\end{equation}
The scalar fields $T^v(\bx,t)$ and $P^v(\bx,t)$ are expanded in series in
the Fourier harmonics $\exp(\ri(\alpha_1n_1x_1+\alpha_2n_2x_2))$; we denote
the coefficients $V^T_{n_1,n_2}(x_3,t)$ and $V^P_{n_1,n_2}(x_3,t)$, respectively.

The boundary conditions (\ref{boum}) imply
\begin{equation}\label{velBC}
T^v|_{x_3=\pm1}=0,\ P^v|_{x_3=\pm1}=0,\
\partial P^v/\partial x_3|_{x_3=\pm1}=0,\ \bM^v|_{x_3=\pm1}=0.
\end{equation}

Computing the r.h.s. of equation (\ref{nst}) is more difficult,
than those of (\ref{mind}) and (\ref{htr}), because it involves $\nabla p$
defined by the condition that the sum is solenoidal. We denote by
$\tilde\bff(\bx,t)$ the r.h.s. of (\ref{nst}) without $\nabla p$,
$$
\tilde\bff={\bf v}\times(\nabla\times{\bf v})+P\tau{\bf v}\times{\bf e}_{\rr}
+P\Delta{\bf v}+PR\theta{\bf e}_z-{\bf b}\times(\nabla\times{\bf b})
$$
and decompose it,
\begin{equation}\label{velrhs}
\tilde\bff(\bx,t)=\bT^f(\bx,t)+\bP^f(\bx,t)+\bM^f(x_3,t)-\nabla p(\bx,t),
\end{equation}
where $\bT^f(\bx,t)$ and $\bP^f(\bx,t)$ are determined by $T^f(\bx,t)$
and $P^f(\bx,t)$ similarly to (\ref{magf1}).
The fields $\tilde\bff(\bx,t)$, $T^f(\bx,t)$ and $P^f(\bx,t)$ are expanded
in series in $\exp(\ri(\alpha_1n_1x_1+\alpha_2n_2x_2))$; we denote
the coefficients $\bF_{n_1,n_2}(x_3,t)$, $G^T_{n_1,n_2}(x_3,t)$ and
$G^P_{n_1,n_2}(x_3,t)$, respectively.

Performing FFT's, multiplication and differentiation, we obtain $\tilde\bff$
in the form of the series (\ref{fdeco}).
By (\ref{magf1}), the coefficients of the components $T^f$ and $\bM^f(x_3,t)$
are determined by the relation
\begin{equation}\label{detcompv}
\renewcommand{\arraystretch}{1.8}
\begin{array}{l}
-(\alpha_1^2n_1^2+\alpha_2^2n_2^2)G^T_{n_1,n_2}=
(\ri \alpha_2n_2 F^1_{n_1,n_2}-\ri \alpha_1n_1 F^2_{n_1,n_2}),\
M_1=F^1_{0,0},\ M_2=F^2_{0,0}.
\end{array}
\end{equation}

Therefore, the Euler scheme for the toroidal and mean-field components is
\begin{equation}\label{eulv}
\renewcommand{\arraystretch}{1.8}
\begin{array}{l}
V^T_{n_1,n_2}(x_3,t_{m+1})=V^T_{n_1,n_2}(x_3,t_m)+
\delta t P_{\cV}(G^T_{n_1,n_2}(x_3,t_m)),\\
M^v_1(x_3,t_{m+1})=M^v_1(x_3,t_m)+\delta t P_{\cV}(F^1_{0,0}(x_3,t_m)),\\
M^v_2(x_3,t_{m+1})=M^v_2(x_3,t_m)+\delta t P_{\cV}(F^2_{0,0}(x_3,t_m)),
\end{array}
\end{equation}
where $G^T_{n_1,n_2}$ are determined by (\ref{detcompv}) and the space $\cV$
is the same as in Example 1.

To determine $G^P_{n_1,n_2}(\bx,t)$, we apply
twice the operator $\nabla\times$ to (\ref{velrhs}) and obtain
\begin{equation}\label{findPv}
\renewcommand{\arraystretch}{1.8}
\begin{array}{l}
(\alpha_1^2n_1^2+\alpha_2^2n_2^2)((\alpha_1^2n_1^2+\alpha_2^2n_2^2)G^P_{n_1,n_2}+\partial^2 G^P_{n_1,n_2}/\partial x_3^2)=\\
(\alpha_1^2n_1^2+\alpha_2^2n_2^2)F^3_{n_1,n_2}+\ri \alpha_1n_1\partial F^1_{n_1,n_2}/\partial x_3+
\ri n_2\partial F^2_{n_1,n_2}/\partial x_3.
\end{array}
\end{equation}
Hence, for the poloidal component of $P_{\cU}\bff_v(\bx,t)$ the computation is
reduced to solving for each pair $(n_1,n_2)$ the equation in
$v=P_{\cV}G^P_{n_1,n_2}$,
$$P_{\cV}(\alpha v+\beta v''-f)=0,\ v\in {\cV},$$
where $\alpha$, $\beta$ and $f$ are determined by (\ref{findPv}) and the
space ${\cV}$ is defined by (\ref{velBC}).
This problem has been considered in Example 2.

\subsection{Numerical results}\label{sec43}

Equations (\ref{nst})-(\ref{htr}) have been integrated in time using
a fourth-order Runge--Kutta scheme
$$
\renewcommand{\arraystretch}{1.8}
\begin{array}{l}
\bd_1=\delta t P_{\cU}\bff(\bu(t_m))\\
\bd_2=\delta t P_{\cU}\bff(\bu(t_m)+\displaystyle{1\over 2}\bd_1)\\
\bd_3=\delta t P_{\cU}\bff(\bu(t_m)+\displaystyle{1\over 2}\bd_2)\\
\bd_4=\delta t P_{\cU}\bff(\bu(t_m)+\bd_3)\\
\bu(t_{m+1})=\bu(t_m)+\displaystyle{1\over6}\bd_1
+\displaystyle{1\over3}\bd_2+\displaystyle{1\over3}\bd_3+
\displaystyle{1\over6}\bd_4
\end{array}
$$
(see \cite{nr}). Here $\bu=(\bv,\theta,\bb)$ and
$f(\bu)=(f_{\bv}(\bu),f_{\theta}(\bu),f_{\bb}(\bu))$ denotes the
r.h.s. of the equations. The projections $P_{\cU}$ are computed
as for the Euler scheme discussed in the previous section.

Results of computations for the parameters
$$P=1,\ R=50000,\ \tau=500,\ \re_{\rr}=(0,1,1),\ P_m=2,$$
the resolution of $32^2$ Fourier harmonics and 16 Chebyshev polynomials, and
the time step $\delta=10^{-4}$ are shown in Fig.~1. The initial condition
is an amagnetic convective steady state of the form of wavy rolls
(the convective attractors in the absence of magnetic field that are
studied in \cite{op24}), and a small magnetic field. Since $P_m=2$ is not far
from the critical value of $P_m$ for the magnetic field generation,
in the nonlinear regime the kinetic energy does not differ much from its value
at $t=0$ and magnetic energy is small relative to the kinetic one.

\vspace*{5mm}
\hspace*{-7mm}\psfig{file=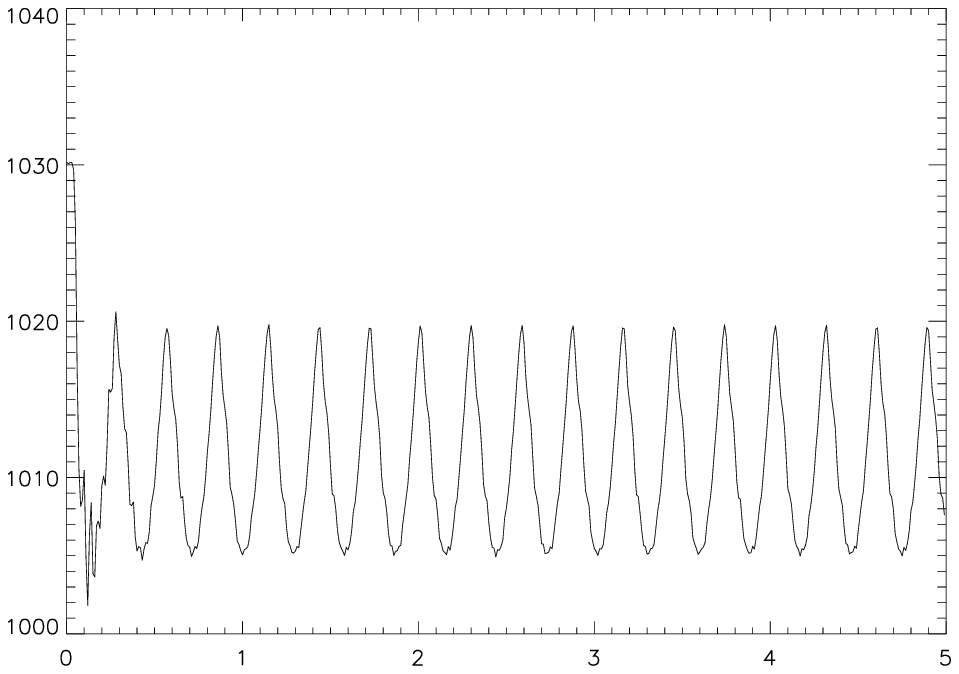,width=8cm}\hspace*{5mm}\psfig{file=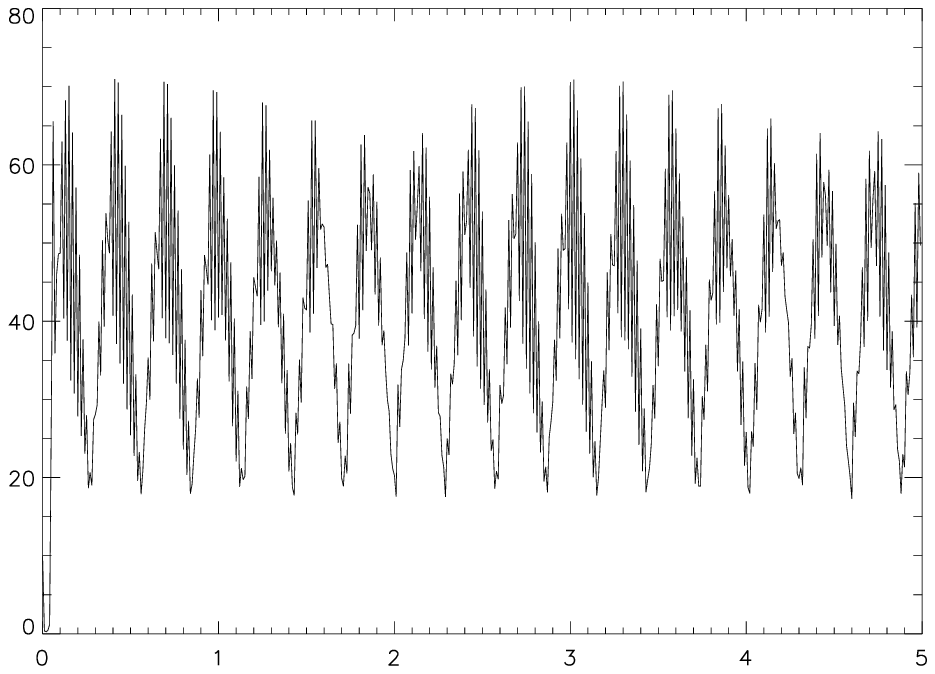,width=8cm}

\vspace*{-5.cm}
\hspace*{-1.1cm}
{\large $E_v$\hspace*{8.2cm}$E_b$}

\vspace*{4.1cm}
\hspace*{5.0cm}{\large $t$}\hspace*{8.3cm}{\large $t$}

\medskip\noindent
Figure 1.
Kinetic (left panel) and magnetic (right panel) energies (the vertical axis) as
functions of time (the horizontal axis).

The aim of this study is investigation of magnetic field generation
by convective attractors in linear and nonlinear regimes, including for $Pm$
substantially larger than the critical values for the onset of the generation.
In order to resolve the emerging small scales in magnetic fields, a higher
space resolution must be employed, which usually results in a reduction
od the time step.
Hence, efficient computing requires the use of implicit-explicit numerical
scheme, which we discuss in the next section.

\subsection{An implicit-explicit scheme}\label{sec44}

In this section we describe an implicit-explicit numerical scheme for
solving the equations of magnetohydrodynamics. Time integration
is performed using the Euler scheme. Generalisations for more complex methods
can be done following the approach of \cite{im1,im2}.

Denote by $f_{\theta}$ the r.h.s. of the heat transfer equation without
the Laplacian
\begin{equation}\label{rhtmin}
f_{\theta}(\bu)=-({\bf v}\cdot\nabla)\theta+v_z.
\end{equation}
For the temperature an implicit-explicit Euler scheme is
$$\theta(\bx,t_{m+1})-\Delta\theta(\bx,t_{m+1})=
\theta(\bx,t_m)+\delta tf_{\theta}(\bu(\bx,t_m).$$
For $\theta$ expanded in a finite sum of products of Fourier exponentials
and polynomials, the Galerkin method yields the equation
\begin{equation}\label{tempim}
P_{\cU}(\theta(\bx,t_{m+1})-\Delta\theta(\bx,t_{m+1}))=
P_{\cU}(\theta(\bx,t_m)+\delta tf_{\theta}(\bu(\bx,t_m)),\
\theta(\bx,t_{m+1})\in\cU,
\end{equation}
where $\cU$ is the space of functions satisfying the boundary condition.

Let $\theta(\bx,t)$ be expanded in series (\ref{temp1}), where we denote
the coefficients $\Theta_{n_1,n_2}(x_3,t)$, and $f_{\theta}(\bx,t)$
in series (\ref{ff1}) with the
coefficients $F_{n_1,n_2}(x_3,t)$. Equation (\ref{tempim}) splits into
$(2N_1+1)\times(2N_2+1)$ equations for the pairs $(n_1,n_2)$,
$$P_{\cV}((1-\delta t(\alpha_1^2n_1^2+\alpha_2^2n_2^2))v-\delta tv''-f)=0,\ v\in {\cV},$$
where $v=\Theta_{n_1,n_2}(x_3,t_{m+1})$, $f=\Theta_{n_1,n_2}(x_3,t_m)+\delta tF_{n_1,n_2}(x_3,t_m)$ and
${\cV}$ is the same as in Examples 1 and 2.

For the magnetic field, we denote by $\bff_b$ the r.h.s. of equation (\ref{mind})
without the Laplacian,
$$\bff_b=\nabla\times({\bf v}\times{\bf b})$$
and expand it in the sum of the toroidal, poloidal and mean-field components
(\ref{magf0}),(\ref{magf1}). The toroidal and poloidal components,
$T^b$ and $P^b$, and $\bff_b$ are expanded as series in $\exp(\ri(\alpha_1n_1x_1+\alpha_2n_2x_2))$
with coefficients $B^T_{n_1,n_2}(x_3,t)$, $B^P_{n_1,n_2}(x_3,t)$ and
$\bF^T_{n_1,n_2}(x_3,t)$, respectively.

By the arguments applied for deriving the implicit-explicit scheme
for the temperature, the numerical scheme for the magnetic field coefficients
takes the form
\begin{equation}\label{umagf}
P_{\cV}((1-\delta t\eta(\alpha_1^2n_1^2+\alpha_2^2n_2^2))v-\delta t\eta v''-f)=0,\ v\in {\cV},
\end{equation}
where
$$v=B^T_{n_1,n_2}(x_3,t_{m+1}),\ f= B^T_{n_1,n_2}(x_3,t_m)+\delta tG^T_{n_1,n_2}(x_3,t_m)$$
for the toroidal coefficients,
$$v=B^P_{n_1,n_2}(x_3,t_{m+1}),\ f= B^P_{n_1,n_2}(x_3,t_m)+\delta tG^P_{n_1,n_2}(x_3,t_m)$$
for the poloidal coefficients, and
$$v=M^j_{n_1,n_2}(x_3,t_{m+1}),\ f= M^j_{n_1,n_2}(x_3,t_m)+\delta tF^j_{0,0}(x_3,t_m),\ j=1,2$$
for the mean field. The coefficients $G^T_{n_1,n_2}$ and $G^P_{n_1,n_2}$ are
determined by (\ref{detcomp}).
The subspaces $\cV$ and $\cV^{\perp}$ are described in section \ref{sec422}.
Equation (\ref{umagf}) is solved employing the same main step as in Example 2,
and two other steps as described in the end of section \ref{sec2}.

For the flow velocity, we denote by $\tilde\bff$ the r.h.s. of the equation
(\ref{nst}) with the pressure and the Laplacian omitted,
$$\tilde\bff={\bf v}\times(\nabla\times{\bf v})+P\tau{\bf v}\times{\bf e}_{\rr}
+PR\theta{\bf e}_z-{\bf b}\times(\nabla\times{\bf b}),$$
and expand it in series (\ref{velrhs}). We decompose the velocity and expand it
in series as discussed in section \ref{sec423}. The numerical scheme
for the toroidal coefficients and the two components of the mean field is
$$P_{\cV}((1-\delta tP(\alpha_1^2n_1^2+\alpha_2^2n_2^2))v-\delta tP v''-f)=0,$$
where for the toroidal coefficients
$$v=V^T_{n_1,n_2}(x_3,t_{m+1}),\ f= V^T_{n_1,n_2}(x_3,t_m)+\delta tG^T_{n_1,n_2}(x_3,t_m)$$
and for the mean field
$$v=M^j_{n_1,n_2}(x_3,t_{m+1}),\ f= M^j_{n_1,n_2}(x_3,t_m)+\delta tF^j_{0,0}(x_3,t_m),\ j=1,2.$$
The coefficients $G^T_{n_1,n_2}$ are determined by (\ref{detcomp}).
The subspaces ${\cV}$ are the same as in Examples 1 and 2.

For the poloidal component, proceeding as in section \ref{sec423}, we find
that $v=V^P_{n_1,n_2}(x_3,t_{m+1})$ can be determined from the equation
$$P_{\cV}((\alpha_1^2n_1^2+\alpha_2^2n_2^2)(1-\delta tP(\alpha_1^2n_1^2+\alpha_2^2n_2^2))v
+(-2\delta tP (\alpha_1^2n_1^2+\alpha_2^2n_2^2)-1) v''+\delta tPv''''-f)=0,\ v\in{\cV}$$
where $f= V^P_{n_1,n_2}(x_3,t_m)+\delta tG^P_{n_1,n_2}(x_3,t_m)$ and
$G^P_{n_1,n_2}$ can be found using (\ref{detcomp}).
The choose the same subspace ${\cV}$ as in Examples 1 and 2. The problem
(the main step) can be solved by an algorithm
similar to that employed to solve the problem in Example 2.

\section{Conclusion}\label{con}

We propose a new algorithm for numerical solution of evolutinary
equations by the Galerkin method. After discretisation in space and time
is implemented, we need to solve equations of the form $P_{\cV}(\cA v-f)=0$,
where $v\in\cV$, at each time step. The equations are solved in three steps:
a preliminary step, which is performed only once; the main step, solving
$P_{\cV}(\cA w-f)=0$, where $w\in\cW$; and computing the correction $v-w$.
With minor modifications the algorithm is also applicable for solving
evolutionary equations by the Petrov--Galerkin method, see remark \ref{rem1}.

The algorithm has been employed for a numerical study of magnetic field
generation by convective flows in a plane layer. The system under
investigation is comprised of the Navier--Stokes equation, the heat transfer
equation and the magnetic induction equation. Representing the solution as
a sum of products of Fourier exponentials in $x_1$ and $x_2$ and Chebyshev
polynomials in $x_3$, we derive a numerical scheme which involves
a set of one-dimensional equations to be solved at each time step.
The method was also used to study convection in a plane
layer rotating about an inclined axis \cite{op24}.

We have presented a detailed derivation of the numerical scheme for the explicit
Euler time stepping. Since explicit multistep schemes also involve
computation of the projections $P_{\cV}f(\bx,t)$ (see, e.g., a discussion of the
Runge--Kutta scheme in section \ref{sec43}), the algorithm can be applied
in the context of other such schemes.

In section \ref{sec44} we show that the algorithm can be used for time stepping
by implicit-explicit methods; we discuss the Euler method as an example.
Although implicit and implicit-explicit methods are more difficult to implement,
they can be more efficient due to the use of larger time steps.
We do not discuss here how the algorithm can be applied to implicit-explicit
multistep methods, this can be done following \cite{asc97,cal01,im1,im2}.

\bigskip
{\large\bf Acknowledgments}
The project was financed by the grant
{\fontfamily{cmr}\selectfont\textnumero}\,22-17-00114 of the Russian
Science Foundation, https://rscf.ru/project/22-17-00114/.

\end{document}